\documentclass[12pt,reqno]{amsart}
\usepackage{amsmath}
\usepackage{amsfonts}
\usepackage{amssymb}
\usepackage{amsthm}
\usepackage{amscd}
\usepackage{a4wide}
\usepackage{enumerate}
\usepackage{dsfont} 

\usepackage{graphics}
\usepackage{eucal}
\usepackage{mathrsfs}

\usepackage{color}
\usepackage[pdftex,bookmarks,colorlinks,breaklinks]{hyperref}  

\hypersetup{linkcolor=red,citecolor=blue,filecolor=dullmagenta,urlcolor=darkblue} 

\newtheorem {theorem}    {Theorem}[section]

\newtheorem {corollary}  [theorem]    {Corollary}

\theoremstyle{definition}
\newtheorem {definition} [theorem]    {Definition}
\newtheorem {example}    [theorem]    {Example}

\newcounter{AbcT}

\numberwithin{equation}{section}

\newcommand {\N} {{\mathbb N}}
\newcommand {\R} {{\mathbb R}}

\newcommand {\Z} {{\mathbb Z}}

\newcommand{\IGNORE}[1]{}

\renewcommand{\limsup}{\varlimsup}

\newcommand\degree{\operatorname{deg}}

\newcommand{\bbZ}{\mathbb{Z}}

\newcommand{\bq}{\mathbf{q}}

\newcommand{\height}{\operatorname{H}}
\newcommand{\theight}{\widetilde{\operatorname{H}}}

\begin{document}
\title{Remarks on Diophantine approximation in function fields}

\author{Arijit Ganguly}
\author{Anish Ghosh}

\address{School of Mathematics, Tata Institute of Fundamental Research, Mumbai, 400005, India}
\email{arimath@math.tifr.res.in, ghosh@math.tifr.res.in}

\date{\today}
\thanks{Ghosh is supported by a UGC grant}

\begin{abstract}
We study some problems in metric Diophantine approximation over local fields of positive characteristic. 

\end{abstract}

\maketitle

\section{Introduction}\label{section:intro}

This short paper is motivated by a recent work of Kristensen, Pedersen and Weiss \cite{KPW}. In it, the authors consider 
certain interesting Diophantine problems arising from K. R. Yu's \cite{Yu} generalisation to higher dimensions, of 
Mahler's \cite{Mah-clas} classification of complex numbers according to their Diophantine properties. In this paper, which can be 
thought of as a sequel to our earlier work \cite{GG} and the paper \cite{G-pos} of the second named author, we consider 
function field analogues of their  results. We recall Mahler's classification briefly for the reader, see \cite{Bug} for 
details. For $k \in \N$ and $\alpha \in \R$ consider the Diophantine exponent

\begin{equation}\label{Mah}
\omega_{k}(x) = \sup\{\omega \in  \R~:~|P(x)| \leq \height(P)^{-\omega}
\text{ for infinitely many irreducible } P \in Z[X], \deg(P) \leq k\}.
\end{equation}

\noindent Here $\height(P)$ denotes the na\"{i}ve height of the polynomial $P$ which is defined in \S 2. Let 

\begin{equation}\label{Mah2}
\omega(x) = \limsup_{k \to \infty} \frac{\omega_{k}(x)}{k}
\end{equation}

Then, according to Mahler, there are four classes of real numbers:
\begin{itemize}
\item $x$ is an $A$-number if $\omega(x) = 0$.
\item $x$ is an $S$-number if $0< \omega(x) < \infty$.
\item $x$ is a $T$-number if $\omega(x)= \infty$, but $\omega_{k}(x)<\infty$ for all $k$.
\item $x$ is a $U$-number if $\omega(x)= \infty$ and $\omega_{k}(x)= \infty$ for all $k$ large enough.

\end{itemize}

Mahler proved that almost all real numbers are $S$-numbers and conjectured \cite{Mah-S} that in fact, $\omega_k(x) = k$ for 
almost every $x \in \R$. This conjecture was proven by Sprindh\v{z}uk \cite{Spr1, Spr2}. It is equivalent to show that almost every 
point on the curve $(x, x^2, \dots, x^n)$ is not very well approximable, i.e. for almost every $x \in \R$, for every $\varepsilon > 0$ there are 
at most finitely many $\bq = (q_1, q_2, \dots, q_n) \in \Z^n$ such that 
$$ |q_1x + q_2x^2 + \dots +q_nx^n +p| \leq \|\bq\|^{-n(1+\varepsilon)}\text{ for some }p\in \bbZ,$$ where $\|\bq\|$ denotes the supremum norm of $\bq$.
\noindent Sprindh\v{z}uk's proof of Mahler's conjecture initiated the subject of ``metric Diophantine 
approximation with dependent quantities" which has seen dramatic progress since then and received an important boost 
with the proof of the Baker-Sprindh\v{z}uk conjectures by D. Kleinbock and G. A. Margulis \cite{KM} using methods from homogeneous dynamics. 

The setting of this present paper involves Diophantine approximation over local fields of positive 
characteristic, i.e. Laurent series over a finite base field. The notation and pertinent definitions are 
introduced in the next section. Diophantine approximation in this setting has been extensively studied beginning 
with Mahler \cite{Mahler, Mahler-L}. Mahler's classification was undertaken in positive characteristic 
by Bundschuh \cite{Bund}. We refer the reader also to \cite{Dub} and to \cite{Ooto} for recent developments 
in analogues of Mahler's and Koksma's classification in positive characteristic. As regards Mahler's conjectures and their 
analogues, Sprindh\v{z}uk  formulated and settled the function field analogue of Mahler's conjectures. The function field 
analogue of the Baker-Sprindh\v{z}uk conjectures were settled by the second named author in \cite{G-pos}. 

In \cite{Yu}, K. R. Yu developed a higher dimensional analogue of Mahler's classification. As far as we are 
aware, this higher dimensional classification has not been investigated in the setting of function 
fields although it seems a natural and interesting problem. In \cite{KPW}, Kristensen, Pedersen and Weiss use Yu's 
classification as a starting point to study various problems arising in metric Diophantine approximation with dependent 
quantities. For instance, they introduce notions of $k$-very well approximable and $k$-Dirichlet improvable vectors and show that 
they are contained in a set of zero Lebesgue measure. In fact their results, as well as ours, are more general; details are in subsequent sections.  In this paper, we develop function field analogues of some of their results.  

We now briefly introduce the setting of Diophantine approximation in function fields. The subsequent 
sections are devoted, respectively, to the notion of $k$ Dirichlet improvable vectors 
and $k$ very well approximable vectors. Dirichlet improvable vectors and very well approximable 
vectors have been studied in the function field setting in a previous paper of the authors \cite{GG} and in a paper \cite{G-pos} of the second named author respectively. Following the ideas in \cite{KPW}, we explain how one can use these previously established results to study the refined notion of $k$-Dirichlet improvable and $k$-very well approximable vectors. We also use the occasion to record a result which characterises Dirichlet improvable Laurent series, which shows a striking contrast with the case of real numbers. Finally, we should mention that \cite{KPW} also has a nice section devoted to badly approximable vectors which is based on the beautiful work of Beresnevich \cite{Ber} on badly approximable vectors on manifolds. Since there are as yet no function field analogues of these results, we have nothing to say about it.

\subsection{The set up}

\noindent Let $p$ be a prime and $q:= p^r$, where $r\in \mathbb{N}$, and consider the function 
field $\mathbb{F}_{q}(T)$. We define a function $|\cdot|: \mathbb{F}_{q}(T) \longrightarrow \mathbb{R}_{\geq 0}$ as follows. 
\[ |0|:= 0\,\,  \text{ and} \,\, \left|\frac{P}{Q}\right|:= e^{\displaystyle \degree P- \degree Q}
\text{ \,\,\,for all nonzero } P, Q\in \mathbb{F}_{q}[T]\,.\] 
Clearly $|\cdot|$ is a nontrivial,  non-archimedian and discrete absolute value  
in $\mathbb{F}_{q}(T)$. This absolute value gives rise to a metric on $\mathbb{F}_{q}(T)$. \\

The completion field of $\mathbb{F}_{q}(T)$ is $\mathbb{F}_{q}((T^{-1}))$, i.e. the field of Laurent series 
over $\mathbb{F}_{q}$. The absolute value of $\mathbb{F}_{q}((T^{-1}))$, which we again denote by $|\cdot |$, is given as follows. 
Let $a \in \mathbb{F}_{q}((T^{-1}))$. For $a=0$, 
define $|a|=0$. If $a \neq 0$, then we can write 
$$a=\displaystyle \sum_{k\leq k_{0}} a_k T^{k}\,\,\mbox{where}\,\,\,\,k_0 \in \mathbb{Z},\,a_k\in \mathbb{F}_{q}\,\,\mbox{and}\,\, a_{k_0}
\neq 0\,. $$
\noindent We define $k_0$ as the \textit{degree} of $a$, which will be denoted by $\degree a$,  and     
$|a|:= e^{\degree a}$. This clearly extends the absolute 
value $|\cdot|$ of $\mathbb{F}_{q}(T)$ to $\mathbb{F}_{q}((T^{-1}))$ and moreover, 
the extension remains non-archimedian and discrete. Let $\Lambda$ and $F$ 
denote $\mathbb{F}_{q}[T]$ and $\mathbb{F}_{q}((T^{-1}))$ respectively from now on. It is obvious that 
$\Lambda$ is discrete in  $F$. For any $d\in \mathbb{N}$, $F^d$ is throughout assumed to be equipped 
with the supremum norm which is defined as follows
\[||\mathbf{x}||:= \displaystyle \max_{1\leq i\leq n} |x_i|\text{ \,\,for all \,} \mathbf{x}=(x_1,x_2,...,x_d)\in F^{d}\,,\]

\noindent and with the topology induced by this norm. Clearly $\Lambda^n$ is discrete in $F^n$. Since the topology on $F^n$ considered here 
 is the usual product topology on $F^n$, it follows that  $F^n$ is locally compact as $F$ is locally compact. Let $\lambda$ be the Haar measure on $F^d$ which takes the value 1 on the closed unit ball $||\mathbf{x}||=1$.\\ 
 
Diophantine approximation in the positive characteristic setting consists of approximating elements 
in $F$ by `rational' elements, i.e. those from $\mathbb{F}_{q}(T)$. As we have mentioned before, this 
subject has been extensively studied, beginning with Mahler who developed Minkowski's geometry of numbers in 
function fields and continuing with Sprind\v{z}uk who, in addition to proving the analogue of Mahler's conjectures, 
also proved some transference principles in the function field setting. The subject has also received 
considerable attention in recent times, we refer the reader to \cite{deM, Las1, Las2} and other works.

 \section{$k$-Dirichlet improvability}Throughout we let $N=\binom {k+d}{d}-1$ be the number of nonconstant monomials in $\Lambda[X_1,...,X_d]$ of total degree $\leq k$, where $k\in \mathbb{N}$, and for $P\in \Lambda[X_1,...,X_d]$, let $\height(P)$  and $\theight(P)$ be the maximum absolute 
 value among all the coefficients of $P$ and the maximum absolute value of the
 coefficients of the non-constant terms of $P$ respectively. We first prove the following analogue of the Dirichlet's theorem in this set up.
\begin{theorem}\label{Dirichlet}
For given $\mathbf{x}=(x_1,...,x_d)\in F^d$, there exists $c(\mathbf{x})$ such that for all $m\in \mathbb{N}$, one has $P\in \Lambda[X_1,...,X_d]$ of 
total degree $\leq k$ and $\height(P)\leq e^m$ satisfying 
\begin{equation}
 |P(\mathbf{x})|\,<\,\frac{c(\mathbf{x})}{e^{mN}}\,.
\end{equation}We can choose $c(\mathbf{x})=1$ if the condition $\height(P)\leq e^m$ is replaced by $\theight(P)\leq e^m$. 
\end{theorem}
The above theorem indeed almost follows from Theorem 2.1 of \cite{GG}. In fact, it is just a restatement of \cite[Theorem 2.1]{GG} if we consider $\theight(P)$. To prove the theorem for $\height(P)$, we shall actually prove the the following which yields Theorem \ref{Dirichlet} at once.

\begin{theorem}\label{thm:D}
Assume $n\in \mathbb{N}, \mathbf{y}:=(y_1,y_2,...,y_n)\in F^n$ and $k:=\max \{0,\deg y_1,\cdots, \deg y_n\}$. For any $m\in \mathbb{N}$
there exist $(q_1,\cdots,q_n,p)\in \Lambda^{n+1} \setminus \{0\}$ 
such that
\begin{equation}\label{eqn:D}|y_1q_1+ y_2q_2+\cdot\cdot\cdot+y_nq_n-p|\,\textless \, \frac{e^{nk}}{e^{nm}}\text{ \,\,and } \displaystyle 
|p|, \max_{1\leq j\leq n} |q_j|\leq e^m\,.\end{equation}
\end{theorem}
\emph{Proof}: If $k=0$ then the conclusion follows immediately from  \cite[Theorem 2.1]{GG}. Suppose now $k>0$. Consider $m>k$. By 
\cite[Theorem 2.1]{GG}, we have $(q_1,\cdots,q_n)\in \Lambda^n \setminus \{0\}$ and $p\in \Lambda$ 
such that
\[|y_1q_1+ y_2q_2+\cdot\cdot\cdot+y_nq_n-p|\,\textless \, \frac{1}{e^{n(m-k)}}\text{ \,\,and } \displaystyle 
\max_{1\leq j\leq n} |q_j|\leq e^{m-k}\,.\] It is easy to see that $|p|\leq e^m$ as $\deg (y_1q_1+ y_2q_2+\cdot\cdot\cdot+y_nq_n)\leq (m-k)+k=m$. So 
(\ref{eqn:D}) holds in this case. When $m=k$, take $q_1=1, q_2=\cdots=q_n=0$ and $p$ to be the polynomial part of $q_1y_1=y_1$ which has degree $\leq k=m$. 
Thus (\ref{eqn:D}) is obvious. We are now left with the case $m<k$. Choose $q_1=q_2=\cdots=q_n=0$ and $p=1$ so that (\ref{eqn:D}) holds again. \,\,\, $\Box$

\begin{corollary}
 For given $\mathbf{x}=(x_1,...,x_d)\in F^d$, there exists $c(\mathbf{x})$ such that 
 \[|P(\mathbf{x})|\,<\,\frac{c(\mathbf{x})}{\height(P)^N}\] holds for infinitely many $P\in \Lambda[X_1,...,X_d]$ of 
total degree $\leq k$.
\end{corollary}

A similar statement holds for $\theight(P)$ and it is easily checked that $c(\mathbf{x})$ can be taken to be $1$ in that case.

The notion of \emph{singular} vectors was introduced by Khintchine \cite{Kh, Cas-book} and the more general notion of \emph{Dirichlet improvable} numbers and vectors was introduced by Davenport and Schmidt \cite{DS1, DS2}. Kristensen, Pedersen and Weiss \cite{KPW}  introduced natural extensions of these  notions, calling them $(k,\varepsilon)$-Dirichlet improvable vectors and $k$-singular vectors. We now introduce the function field analogues. Let $0<\varepsilon\leq\frac{1}{e}$. A point 
$\mathbf{x}\in F^d$ is said to be $(k,\varepsilon)-Dirichlet\,improvable$ if there exists $m_0\in \mathbb{N}$ such that for 
every $m\geq m_0$, the following system of inequations admits a nonzero solution $P\in \Lambda[X_1,...,X_d]$ with total degree $\leq k$
\[\left \{\begin{array}{rcl}
|P(\mathbf{x})|< \displaystyle \frac{\varepsilon}{e^{mN}}\\ \theight(P)< \varepsilon e^m
\end{array}\right.\,. \] When $k=1$, this coincides with the usual Dirichlet improvable vectors in $F^d$, studied in \cite{GG}. Furthermore, 
one has the following observation in the case $d=1$. 
\begin{theorem}\label{thrm:DI are rationals}
 For any $\alpha \in F$, $\alpha$ is Dirichlet $\varepsilon$-improvable for some $0 < \varepsilon \leq \frac{1}{e}$ if and only if it is a rational function.\end{theorem}
 \emph{Proof} : \noindent Let $\alpha \in F$. Assume that the Dirichlet's theorem can be $\varepsilon$-improved for $\alpha$, where 
$0\,\textless\, \varepsilon\leq \frac{1}{e}$. So there exists $N_0\in \mathbb{N}$ such that 
for any $n\geq N_0$, one has $p_n,q_n \in \Lambda$ with 
\begin{equation}\label{eqn:DI1}
 \left\{ \begin{array}{rcl} |q_n\alpha- p_n|\,\textless\, \frac{\varepsilon}{e^n} \\ 1\leq |q_n|\,\textless \,\varepsilon e^{n}\,. 
 \end{array} \right.
\end{equation}
As a consequence, we have 
\begin{equation}\label{eqn:DI2}
 \forall~n\in \mathbb{N}\text{ with }n\geq N_0, \left\{ \begin{array}{rcl} |q_n\alpha- p_n|\,\textless\, \frac{1}{e^{n+1}} \\ 1\leq |q_n|\,\textless \, e^{n-1} \end{array} \right.
\end{equation}
 as $\varepsilon \leq \frac{1}{e}$. Without any loss in generality, we assume always that the pairs $p_n, q_n$ are relatively 
 prime and $q_n$ is monic. 
 Now observe that, for each such $n$,
 \[q_np_{n+1}-q_{n+1}p_n= q_{n+1}(q_n\alpha- p_n)- q_n (q_{n+1}\alpha- p_{n+1})\,.\]
 Taking the absolute value of the above and using the ultrametric property,
\begin{equation}\label{the equation}
 |q_np_{n+1}-q_{n+1}p_n|\leq \max \{|q_{n+1}||q_n\alpha- p_n|,|q_n||q_{n+1}\alpha- p_{n+1}|\}\,.
\end{equation} In view of (\ref{eqn:DI2}), it is obvious that, $\forall~n\in \mathbb{N}\text{ with }n\geq N_0$,
\[|q_{n+1}||q_n\alpha- p_n|\,\textless\, \frac{e^n}{e^{n+1}}=\frac{1}{e}\,\textless\,1\] and 
\[|q_n||q_{n+1}\alpha- p_{n+1}|\,\textless\, \frac{e^{n-1}}{e^{n+2}}=\frac{1}{e^3}\,\textless \,1\,.\]
 This shows that the RHS of (\ref{the equation}) is $\textless\,1$ and hence 
 \[q_np_{n+1}-q_{n+1}p_n=0 \text{ i.e. } \frac{p_n}{q_n}=\frac{p_{n+1}}{q_{n+1}}\text{ for all }n\geq N_0\,.\]
Since all $q_n$ are monic and the pairs pairs $p_n, q_n$ are relatively 
 prime, it follows that $p_n=p_{N_0}$ and $q_n=q_{N_0}$ for all $n\geq N_0$. Therefore, from (\ref{eqn:DI2}), one obtains 
 \[|q_{N_0}\alpha -p_{N_0}|\,\textless\,\frac{1}{e^{n+1}} \text{ for all }n\geq N_0\,,\]
and this makes $\alpha = \displaystyle \frac{p_{N_0}}{q_{N_0}}\in \mathbb{F}_q (T)$. Thus if at all any improvement of Dirichlet's theorem
is possible for $\alpha$ then it has to be a rational function. On the other hand, any rational function is trivially
Dirichlet $\varepsilon$-improvable for any $\varepsilon\in (0,\frac{1}{e}]$. \,\,$\Box$ 

We note that in the case of real numbers, the situation is quite different, 
namely there do exist irrational real numbers which are Dirichlet impriovable. In fact,
it was shown by Davenport and Schmidt \cite{DS1} that an irrational number is Dirichlet improvable if and only if it is badly approximable. But for function field with 
positive characteristic, the sets consisting of badly approximable and Dirichlet Laurent series respectively are disjoint.

The set of all $(k,\varepsilon)$-Dirichlet improvable vectors will be denoted by $DI(k,\varepsilon )$. We  call a  
vector $\mathbf{x}\in F^d$ $k-singular$ if it is  $(k,\varepsilon)$-Dirichlet improvable for 
every $0<\varepsilon \leq \frac{1}{e}$. As an obvious corollary of \cite[Theorem 3.7]{GG}, we have 
$\lambda (DI(1,\varepsilon ))=0$ for $\varepsilon\ll1$. 

In this section, we aim to establish that for any $k\in \mathbb{N}$, one has 
$\varepsilon_0= \varepsilon_0 (d,k)>0$ such that $\lambda (DI(k,\varepsilon ))=0$ whenever $\varepsilon <\varepsilon_0$. In fact, we will be proving 
a much more general result. For proceeding towards that, we need a few definitions. \\

Let $\mu$ be a Radon measure
on $F^d$ and $U\subseteq F^d$ be  open with $\mu(U) > 0$. 
The measure $\mu$ is said to be $Federer$ on $U$ if there exists $D>0$ such that for every ball $B$ with center 
 	in $\text{supp }(\mu)$ such that $3B\subseteq U$, one has \[\frac{\mu(3B)}{\mu(B)}\leq D\,.\] Consider a continuous function
 	$f: U\longrightarrow F$. For any $B\subseteq U$, we set
$$ ||f||_{\mu,B} := \displaystyle \sup_{\mathbf{x}\in B\cap \text{ supp }(\mu)} |f(\mathbf{x})|.$$

\begin{definition}\label{defn:C,alpha}
	For $C,\alpha \textgreater\,0$, the function $f$ is said to be $(C,\alpha)-good$ on $U$ with respect to $\mu$ if for every ball $B\subseteq U$
	with center in $\text{supp }(\mu)$, one has
	\[\mu(\{\mathbf{x}\in B: |f(\mathbf{x})|\,\textless \varepsilon\})\leq C\left(\frac{\varepsilon}{||f||_{\mu,B}}\right)^{\alpha} \mu(B)\,.\]
\end{definition}
 We say a  map $\mathbf{f}=(f_1,f_2,...,f_n)$ from $U$ to $F^n$, where $n\in \mathbb{N}$,
 is $(C,\alpha)-good$ on $U$ with respect to $\mu$, or simply $(\mathbf{f},\mu)$ is 
 $(C,\alpha)-good$ on $U$,  if every $F$-linear combination of $1,f_1,...,f_n$ 
 is $(C,\alpha)-good$ on $U$ with respect to $\mu$. \\
 
 \begin{definition} Let $\mathbf{f}=(f_1,f_2,...,f_n)$ be a map from $U$ to $F^n$, where $n\in \mathbb{N}$. 
 	We say that $\mathbf{f}$ is \emph{nonplanar} with respect to $\mu$ or $(\mathbf{f},\mu)$ is \emph{nonplanar} if for any ball $B\subseteq U$ with center in $\text{supp }(\mu)$, 
 	the restrictions of the functions  $1,f_1,...,f_n$ 
 	on $B\cap \text{ supp }(\mu)$ are linearly independent. 
 \end{definition}
 
The notion of nonplanar maps was introduced in \cite{KT} and generalizes the notion of nondegenerate maps introduced in \cite{KM}.
 
For $m\in \mathbb{N}$ and a ball $B=B(\mathbf{x};r)\subseteq F^d$, where $\mathbf{x}\in F^d$ and $r\,\textgreater\,0$, we shall use the notation
 $3^mB$ to denote the ball $B(\mathbf{x};3^mr)$. 
 
 \begin{definition}
 	We call $\mu$, $k$-friendly on $U$ if it is Federer on $U$ and the 
 	function 
 	\begin{equation}\label{func}\mathbf{f}:U \longrightarrow F^N
 	\text{ given by } (x_1,x_2,\cdots,x_d)\mapsto (x_1,x_2,\cdots,x_{d-1}x_d^{k-1}, x_d ^k),\end{equation} 
 	i.e. $\mathbf{f}$ maps $(x_1,x_2,\cdots,x_d)$ to $N$-distinct monomials in $d$-variables of total degree $\leq k$, is $(C,\alpha)$-good for some $C,\alpha >0$ and nonplanar with respect to $\mu$.   
 \end{definition}
 
Our definition of $k$-friendly measures is slightly weaker than that in \cite{KPW} leading to potentially stronger results. We refer the reader to \S 10.5 of \cite{KT} where these definitions are compared. The notion of friendly measures was originally introduced in \cite{KLW}.
 
 \begin{example}
 It is not dificult to see that the  prefixed Haar measure $\lambda$  on $F^d$ is Federer on any given open $U\subseteq F^d$ with $D=e^2$ and 
 nonplanar. By \cite[Lemma 2.4]{KT}, $\mathbf{f}$ is 
 $(C,\frac{1}{dk})$-good, where $C$ depends on $d$ and $k$ only. Thus $\lambda$ is $k$-friendly on any open $U\subseteq F^d$. 
 \end{example}
 With these definitions, we now discuss the main theorem of this section.
 \begin{theorem}\label{improv}
  Suppose $\mu$ is a $k$-friendly measure  on $U\subseteq F^d$. Then there exist $\varepsilon_0=\varepsilon_0(d,\mu)$ such that  
  $\mu (DI(k,\varepsilon ))=0$ whenever 
  $\varepsilon <\varepsilon_0$. Thus in particular, the set of all $k$-singular vectors is $\mu$-null.
 \end{theorem}
 \emph{Proof}: The key observation is that if $\mathbf{x}\in U$ is $(k,\varepsilon)$-Dirichlet improvable, 
 then $\mathbf{f}(\mathbf{x})$, where $\mathbf{f}$ is defined in (\ref{func}), is a Dirichlet $\varepsilon$-improvable vector in $F^{N}$ 
 in the sense of \cite{GG}. To see this, let us write the $N$ monomials in the variables $X_1, \cdots, X_d$ having total degree $\leq k$ as $$M_1(X_1, \cdots, X_d), 
 \cdots ,M_N(X_1, \cdots, X_d)$$ so that a polynomial $P\in \Lambda[X_1, \cdots, X_d]$ will look like 
 \[a_0+a_1M_1(X_1, \cdots, X_d)+
 \cdots +a_NM_N(X_1, \cdots, X_d),\]
 where $a_0,\cdots,a_N\in \Lambda$. Clearly $\mathbf{f}(\mathbf{x})=(M(\mathbf{x}), \cdots, M_N(\mathbf{x}))\, \forall~\mathbf{x}\in F^d$. From this, one can see that, if 
 $\mathbf{x}\in U$ is $(k,\varepsilon)$-Dirichlet improvable, where $\varepsilon\in (0, 1/e]$, then for all $m\gg1$, there exist $a_0,\cdots,a_N\in \Lambda$ such that 
 \[\left \{\begin{array}{rcl}
|a_0+(a_1, \cdots, a_N)\cdot\mathbf{f}(\mathbf{x})|=|a_0+a_1M_1(\mathbf{x})+
 \cdots +a_NM_N(\mathbf{x})|< \displaystyle \frac{\varepsilon}{e^{mN}}\\ ||(a_1, \cdots, a_N)||< \varepsilon e^m
\end{array}\right.\,. \]This establishes the claim. We therefore apply \cite[Theorem 3.7]{GG} with $\mathbf{f}$ given by (\ref{func}). From the definition 
 of $k$-friendly measure, the hypothesis of \cite[Theorem 3.7]{GG} is already satisfied, whence considering 
 \[\mathcal{T}:= \{(mN,m,m,\cdots,m):m\in \mathbb{N}\},\]
 Theorem \ref{improv} follows.\,\,$\Box$ 
\section{$k$-Very well approximable points}
A point 
$\mathbf{x}\in F^d$ will be called $k-very\,well\,approximable$ if for some $\varepsilon>0$, one has 
\[|P(\mathbf{x})|<\displaystyle \frac{1}{\height(P)^{N+\varepsilon}}\text{ for infinitely many } P\in \Lambda[X_1,...,X_d]\,.\] We aim to prove that this property is exceptional in the sense that the set of all such points in $F^d$ is $\lambda$-null. In fact, we prove this for any $k$-friendly measure. 
\begin{theorem}
	For any $k$-friendly measure $\mu$ on a given $U\subseteq F^d$, The set of all $k$-VWA points in $U$ is null with respect to $\mu$. 
\end{theorem}
The main observation to prove this theorem is the following. For given $\mathbf{x}\in U$, denote $\mathbf{y}=\mathbf{f}(\mathbf{x})$, where $\mathbf{f}$ is given by (\ref{func}). Now for any polynomial $P$ having total degree $\leq k$, we have $P(\mathbf{x})=\mathbf{q}\cdot \mathbf{y}+p$ where $\mathbf{q},p\in \Lambda$, so that $\height(P)\geq ||\mathbf{q}||$. Hence, if $\mathbf{x}\in F^d$ is a $k$-VWA point then $\mathbf{f}(\mathbf{x})$ is a VWA vector in the sense of \cite[Definition 1.1]{G-pos}. Thus the above theorem follows imemdiately from \cite[Theorem 1.6]{G-pos} when $\mu=\lambda$. For the general case, namely to consider more general measures than Haar measure, one has to adopt the modifications considered in \cite{KT}.

\end{document}